\documentclass{article}
\usepackage[russian]{babel}
\usepackage{amssymb, amsmath, euscript}
\begin{document}
\large
\begin{center}

{\bf\LARGE Statistical Investigation of Increments

of Currency Rates Logarithms}

\bigskip

{\it 17th International Conference "Lomonosov - 2010"\ }

{\it A. Sarantsev}

{\it University of Washington, PhD Student}

{\it E-mail: ansa1989@u.washington.edu}

\end{center}

The aim of this report is to study the behavior of some statistical features of currency rates dynamics. We shall also analyze the applicability of the celebrated Samuelson model and its modification with normal inverse Gaussian (NIG) increments of logarithms to currency rates dynamics modeling. We shall always work not with currency rates themselves, but with their logarithms increments.

We took rates of 12 currencies from the website [6]. These are the dollar, the euro, the Australian dollar, the Canadian dollar, the Iceland krona, the Kazakhstan tenge, the Norway krona, the SDR, the Singapore dollar, the UK pound, the Ukranian griven, the Japanese yen. The rates ar taken with respect to Russian rouble for the period from January 1, 2000, to December 31, 2009. One value of each currency rate is taken for every working day.

Each of these rates is analyzed for each of the 20 half-years which comprise this 10-year period. In addition, similar analysis is carried out separately for the period before the financial crisis and the period during this crisis. We conventionally consider August 1, 2008 as the beginning of the crisis.

The simplest model for currency rate is the Samuelson model $X(t) = e^{at + \sigma W(t)}$, where $W$ is the standard Brownian motion. But here we operate with discrete time and hence we deal with the sequence $X_0, \ldots, X_N$, where $x_k = \ln(X_{k+1}/X_k), k = 0, \ldots, N - 1$ is a selection from the normal distribtion $\mathcal N(a; \sigma^2)$.

{\bf 1. Empirical Volatility.} $$s^2 = \frac1{n-1}\sum_{k=0}^{n-1}(x_k - \overline x)^2.$$ This $s^2$ dramatically rose during the crisis (its order substantially increased). One may consider the volatility as an indicator of weight of the financial situation. It is noteworthy that different currencies have very different $s^2$ (with different orders).

This Samuelson model has a drawback, namely that $x_k$ are Gaussian. Practice shows that the tails of these variables have heavier tails. We shall test whether these variables are Gaussian. One could estimate mean and variance and apply the standard Kolmogorov test. But we shall perform other tests.

{\bf 2. Empirical Asymmetry and Excess Coefficients.} The results of this test show a relatively poor conformation with the assumption that $x_t$ are normal. These coefficients show the largest deviation from the theoretical values for the dollar during the pre-crisis period and for the Kazakhstan tenge during the crisis period.

{\bf 3. Pickands Estimate.} (Cf. [5], chapter 3.) This is an estimate of the extreme value index of a selection distribution. This index is equal to $0$ for the normal distribution. These results do not conform with the Samuelson model too, and one can observe the lack of statistical stability. However, this estimate does not show the more irregular behavior during the crisis than during the pre-crisis period.

{\bf 4. Hill Estimate.} (Ñf. [5], chapter 3.) If the extreme value index is positive, the Hill estimate can be also applied. It is recommended in [3] to test heavy tails using this estimate. For the given $x_t$ one can conclude: $x_t$ has a positive extreme value index (within the range $0.3 - 0.8$, it slightly increased during the crisis). Also, $x_t$ has a power tail, i.e. is not normal. In addition, a statistical stability is not observed.

Another drawback of the Samuelson model is that $x_t$ are independent. There is an opinion that this is not always true. How shall we test whether $x_t$ are independent? Let us suggest the following test, which we call an M-test. (See also [2].)

{\bf 5. M-test.} Suppose $(X_n, n \in \mathbb Z)$ is a sequence of independent identically distributed random variables with continuous distribution function. Let us call a point $n \in \mathbb Z$ a {\it local maximum point} if $X_n > X_{n - 1},\ X_n > X_{n + 1}$. Similarly, let us call this point a {\it local minimum point} if $X_n < X_{n + 1},\ X_n < X_{n - 1}$.

It can be easily shown that local maximum points form an infinite set with infinite upper and lower bounds. The same is true for local minimum points. These points alter, i.e. there is a unique point of local minimum between any two adjacent points of local maximum, and vice versa.

Suppose $0$ is a local maximum point. Let $\tau_1$ be the distance from this point to the next local minimum point, let $\tau_2$ be the distance from this point to the next local maximum point, etc. It is straightforward to prove that $\tau = (\tau_n, n \in \mathbb N)$ is a strictly stationary sequence and $m := \mathbf E\tau_1 = 3/2$.

{\it Conjecture.} This sequence is ergodic.

If this is true then this sequence satisfies the Birkhoff Ergodicity Theorem: $$\hat{m}_n := \frac1n\sum_{k=1}^n\tau_k \to m$$ almost surely and in the mean sense. The Monte-Carlo modleing showed that this Law of Large Numbers is indeed true; in addition, the following limit theorem is true: $\xi_n := \sqrt{n}(\hat{m}_n - m) \stackrel[d]{\longrightarrow} Y$, where $\mathbf{E} Y = 0,\ Var Y \approx 0.73$. We have also found some fractiles of $Y$. It is now obvious how to construct an i.i.d. (and continuity of distribution function) test.

Test results do not entirely conform with the Samuelson model. The dollar has the most irregular behavior. It is noteworthy that the results during the crisis do not deviate from $m$ more than the results during the pre-crisis period.

{\bf 6. Self-similarity Parameter Estimate.} The R/S-method is proposed in [3, 4] to estimate the self-similarity Hurst parameter for the given process $(X_n, n \in \mathbb Z_+)$ and to test whether this process has any self-similarity.
Let us make its brief outline. Let $$S^2(n) = \frac1{n-1}\sum\limits_{k=0}^{n-1}(X_k - \overline{X})^2$$ be an empirical variance for the observations $X_0, \ldots, X_{n-1}$ and let $S_k := \sum_{j=0}^{k-1}X_k$. The variable $$R(n) = \max_{k = 0, \ldots, n - 1}\biggl(S_k - \frac kn S_n\biggr) - \min_{k = 0, \ldots, n - 1}\biggl(S_k - \frac kn S_n\biggr)\ \ - $$ is called {\it the range}. It is known that {\it the standardized range}
$$\frac{R(n)}{S(n)} \backsim cn^H.$$ as $n \to \infty$. Hence one can estimate $H$ by taking logarithms of both parts of this approximate equality and applying the method of minimal squares. We have $H = 1/2$ for $X_k := x_k$, where $X_k$ are from the Samuelson model.

The following results were obtained: the calculated estimates $\hat{H}$ contradict the supposition that $H$ is equal to $1/2$. The estimate is substantially higher during the crisis than for the pre-crisis period. (The only exception is the Ukranian griven: $\hat{H} \approx 0.927$ before the crisis, $\hat{H} \approx 0.440$ during the crisis.)

All these results clearly show that currency rates dynamics cannot be satisfactorily described by the Samuelson model.

{\bf 7. Normal Inverse Gaussian Distribution.} Using non-Gaussian (e.g., generalized hyperbolic, GH) distributions for $x_t$ is a modern approach.  Istigecheva (2006, p. 11) claims that these distributions can be fruitful in financial time series modeling. They have attractive properties: they are not necessarily symmetric and their tails are heavier than normal tail.

Normal Inverse Gaussian (NIG) distribution is a particular case of GH-distribution. Its density is
$$f(x; \alpha, \beta, \delta, \mu) := \frac{\alpha\delta}{\pi}\frac{K_1\left(\alpha\delta\sqrt{1 + \left(\frac{x - \mu}{\delta}\right)^2}\right)}{\sqrt{\delta^2 + (x - \mu)^2}}\exp\left(\delta\sqrt{\alpha^2 - \beta^2} + \beta(x - \mu)\right).$$
It depends on four parameters: $\alpha, \beta, \delta, \mu$. Here $K_1$ is a modified Bessel function:
$$K_1(x) := \frac12\int_0^{+\infty}\exp\left(-\frac x2\left(y + \frac1y\right)\right)dy.$$
To estimate these parameters, the maximal likelihood estimate is used in the article [1]. But we shall use the method of moments. Suppose $E, V, S, K$ is an expectation, variance, asymmetry and excess coefficients, then (see [1])
$$E = \mu + \delta\tau,\ V = \frac{\delta^2(1 + \tau^2)}{\zeta},\ S = \frac{3\tau}{\sqrt{\zeta(1 + \tau^2)}},\ K = \frac3{\zeta}\left(1 + 4\frac{\tau^2}{1 + \tau^2}\right).$$
Here we denote $\beta = \zeta\tau/\delta, \alpha = \zeta\sqrt{1 + \tau^2}/\delta$. One only needs to calculate estimates for $E, V, S, K$ and to solve this system of equations.

The results: this model is not suitable for currency rates dynamics, since $\zeta$ has a strong statistical instability.

\bigskip

\centerline{\bf References}

1. E. V. Istigecheva. Estimation of Hyperbolic and Inverse Gaussian Distribution Parameters. // Isvestia Tomskogo Polytechnitcheskogo Universiteta, volume 309, issue 6. Tomsk, Publishing House of Tomsk Polytechnic University, 2006. Pages 11-13. (In Russian)

2. A. A. Sarantsev. On a Generalization of Bernoully and Euler Numbers. // Proceedings of the 10th International Seminar on Discrete Mathematics and its Applications. Publishing House of the Department of Mechanics and Mathematics, Lomonosov Moscow State University. Moscow, 2010. (In Russian)

3. O. I. Shelukhin, A. V. Osin, S. M. Smolskiy. Self-Similarity and Fractals. Applications in Telecommunication. Moscow, Fizmatlit, 2008. (In Russian)

4. A. N. Shiryaev. The Essentials of Stochastic Finance. World Scientific, 2003.

5. L. de Haan, A. Ferreira. Extreme Value Theory. An Introduction. Springer, 2006.

6. www.bankir.ru

\end{document}